\documentclass[preprint,authoryear,12pt]{elsarticle}
\usepackage{amsmath}
\usepackage{amsthm}

\newtheorem{theorem}{Theorem}
\newtheorem{lemma}[theorem]{Lemma}
\newtheorem{remark}[theorem]{Remark}
\biboptions{square}
\journal{Stochastic Processes and their Applications}

\begin{document}

\begin{frontmatter}



\title{Truncated variation, upward truncated variation 
and downward truncated variation of Brownian motion with drift 
- their characteristics and applications}


\author{Rafa\l\  Marcin \L ochowski \corref{cor1}}
\cortext[cor1]{Tel.: +48 22 564 9257; fax: +48 22 564 9257 \\
\textit{Email address:} rlocho@sgh.waw.pl \\
\textit{URL:} http://akson.sgh.waw.pl/\textasciitilde rlocho/indexeng.html}

\address{Department of Mathematics and Mathematical Economics, \\ 
Warsaw School of Economics, \\ Al.~Niepodleg\l o\'{s}ci 164, 02-554 Warszawa (Warsaw), Poland}

\begin{abstract}
In \cite{L2008} we defined truncated variation of Brownian motion with
drift, $W_t = B_t + \mu t, t\geq 0,$ where $(B_t)$ is a standard Brownian
motion. Truncated variation differs from regular variation by neglecting jumps smaller than some fixed $c>0.$ We prove that truncated variation is a random variable with finite moment-generating function for any complex argument. 

We also define two closely related quantities - upward truncated variation and downward truncated variation. 

The defined quantities may have some interpretation in financial mathematics. Exponential moment of upward truncated variation may be interpreted as the maximal possible return from trading a financial asset in the presence of flat commission when the dynamics of the prices of the asset follows a geometric Brownian motion process.

We calculate the Laplace transform with respect to time parameter of the moment-generating functions of the upward and downward truncated variations. 

As an application of the obtained formula we give an exact formula for expected value of upward and downward truncated variations. We give also exact (up to universal constants) estimates of the expected values of the mentioned quantities. 
\end{abstract}


\end{frontmatter}

% \linenumbers


\section{Introduction}

\label{Intr}

Let $(W_t, t \geq 0)$ be a Brownian motion with drift, $W_{t}=B_{t}+\mu t,$
where $(B_t, t \geq 0)$ is a standard Brownian motion.

The well known result of Paul L\'{e}vy (cf. \cite{Levy}) states that for any 
$0 \leq a <b$ and any $p \leq 2$ the $p$-variation of the process $W_t$ on
the interval $[a, b]$ is almost surely infinite: 
\begin{equation*}
\sup_{n}\sup_{a\leq t_{1}<t_{2}<...<t_{n}\leq b}\sum_{i=1}^{n-1}\left|
W_{t_{i+1}}-W_{t_{i}}\right|^p=+\infty
\end{equation*}
and if $a \leq t_{1, k} < t_{2, k} < \ldots < t_{n_k, k} \leq b $ is a descending sequence of partitions of the interval $[a,b]$ such that $\lim_{k
\rightarrow \infty} \max_{1 \leq i \leq n_k -1} (t_{i+1,k} - t_{i,k}) = 0,$
then 
\begin{equation}
\lim_{k \rightarrow \infty} \sum_{i=1}^{n_k-1}\left|
W_{t_{i+1,k}}-W_{t_{i,k}}\right|^2= b-a \text{ a.s.}  \label{levy}
\end{equation}

The further results of this type state that if $n_k \rightarrow \infty$ and 
$\max_{1 \leq i \leq n_k -1} (t_{i+1,k} - t_{i,k}) = o(1/\ln(n_k))$ then
equality (\ref{levy}) also holds (\cite{Dudley}), but if it is not true,
then (\ref{levy}) may not be true as well (\cite{delaVega}).

In 1972 S. J. Taylor proved (\cite{Taylor1972}) that the function $\psi(x) =
x^2/ \ln \max\left\{\ln(1/x), e \right\}$ is a function with the smallest
order around $0$ and such that

\begin{equation*}
\sup_{n}\sup_{a\leq t_{1}<t_{2}<...<t_{n}\leq b}\sum_{i=1}^{n-1} \psi\left(
\left| W_{t_{i+1}}-W_{t_{i}}\right| \right) < +\infty \text{ a.s.}
\end{equation*}

In the paper \cite{L2008} we started to investigate another type of
variation of Brownian paths, which neglects small jumps (smaller than some 
$c>0$) and defined \emph{truncated variation} of $W_t$ at the level $c>0$ on
the interval $\left[ a,b\right] $ as 
\begin{equation*}
TV_{\mu }^{c}\left[ a,b\right] := \sup_{n}\sup_{a\leq t_{1}<t_{2}< \ldots
<t_{n}\leq b} \sum_{i=1}^{n-1} \phi_c \left( \left| W_{t_{i+1}}-W_{t_{i}}
\right| \right),
\end{equation*}
where $\phi_c (x) = \max\left\{x -c, 0\right\}.$ We will prove that the truncated variation is not only finite almost surely, but has finite moment-generating function for any complex number.

\begin{remark}
A. N. Chuprunov pointed to the author that it would be also interesting to
have estimates of quadratic truncated variation, which one may define as 
\begin{equation*}
QTV_{\mu }^{c}\left[ a,b\right] :=\sup_{n}\sup_{a\leq t_{1}< \cdots < t_{n}\leq b}\sum_{i=1}^{n-1}\phi_{c^2} \left( \left| W_{t_{i+1}}-W_{t_{i}}
\right|^2 \right) .
\end{equation*}
\end{remark}

\begin{remark}
Similar concept of truncation (or \emph{shrinking}) of random variables on Hilbert spaces investigated Z. Jurek in series of his papers beginning with \cite{J 1975}, \cite{J 1985}, which now evolved into the theory of s-selfdecomposable distributions (see e.g. \cite{J 2004}). 
\end{remark}

Let us define two quantities closely related to truncated variation - \emph{upward truncated variation} of $W_t$ on the interval $[a,b]$
\begin{equation*}
UTV_{\mu }^{c}\left[ a,b\right] :=\sup_{n}\sup_{a\leq
t_{1}<s_{1}<t_{2}<s_{2}<...<t_{n}<s_{n}\leq b}\sum_{i=1}^{n} \phi_c\left(
W_{s_{i}}-W_{t_{i}}\right)
\end{equation*}
and, analogously, \emph{downward truncated variation}
\begin{equation*}
DTV_{\mu }^{c}\left[ a,b\right] :=\sup_{n}\sup_{a\leq
t_{1}<s_{1}<t_{2}<s_{2}<...<t_{n}<s_{n}\leq b}\sum_{i=1}^{n} \phi_c \left(
W_{t_{i}}-W_{s_{i}}\right).
\end{equation*}

The defined quantities are related in the following way 
\begin{eqnarray}
\max \left\{ UTV_{\mu }^{c}\left[ a,b\right] ,DTV_{\mu }^{c}\left[ a,b\right]
\right\} &\leq &TV_{\mu }^{c}\left[ a,b\right]  \notag \\
&\leq &UTV_{\mu }^{c}\left[ a,b\right] +DTV_{\mu }^{c}\left[ a,b\right] .
\label{rel}
\end{eqnarray}

It is easy to see that the three above defined quantities have the following properties, which we state only for the truncated variation

\begin{enumerate}
\item  Shift invariance property in distributions: for any stopping time $\Delta$ relative to the natural filtration of $(W_t, t \geq 0)$
\begin{equation*}
\mathcal{L} \left( TV_{\mu }^{c}\left[ a,b\right] \right) = \mathcal{L} \left( TV_{\mu }^{c}\left[ a+\Delta ,b+\Delta \right] \right).
\end{equation*}

\item  Superadditivity property: for any numbers $a\leq a_{1}<a_{2}<\cdots
<a_{n}\leq b$ 
\begin{equation*}
TV_{\mu }^{c}\left[ a,b\right] \geq \sum_{i=1}^{n-1}TV_{\mu }^{c}\left[
a_{i},a_{i+1}\right] .
\end{equation*}
\
\end{enumerate}

Upward truncated variation has also some interpretation in financial
mathematics. We will prove that $\exp UTV_{\mu }^{c}\left[ a,b\right] -1$ is the least upper bound for the maximum possible rate 
of return from any trading a single asset on time interval $\left[ a,b\right]$ in the presence of flat commission (proportional
to the value of the transaction) when asset's prices follow the geometric
motion process $\exp \left( W_{t}\right).$

Due to this fact and (\ref{rel}) we will be interested in calculating the
moment-generating function of the variables $UTV_{\mu }^{c}\left[ a,b\right]$ and $DTV_{\mu }^{c}\left[ a,b\right].$

Since the distribution of $DTV_{\mu }^{c}\left[ a,b\right]$ is the same as
the distribution of $UTV_{-\mu }^{c}\left[ a,b\right],$ we will deal with
the moment-generating function of upward truncated variation only.

More precisely, we will find the Laplace transform with respect to time
parameter $T$ of the ... moment-generating function of the variable $UTV_{\mu }^{c}\left[ 0,T\right].$ Let us explain that here we use term "Laplace transform" in a broad sense. 
For a measurable (with respect to the Lebesgue measure $dt$) complex function $f,$ defined on a positive half-line, by the Laplace transform of $f$ we will mean the value of the integral $\int_0^{\infty} e^{\nu t} f(t) dt$ for any complex $\nu,$ for which this integral exists. Similarly, by the moment-generating function of a complex random variable $X$ we will mean the expected value $\mathbf{E}\exp (\lambda X)$ for any complex $\lambda,$ for which this value is well defined.

As an application of the obtained formula we will give an exact formula for
expected value of upward and downward truncated variations. We give also
exact (up to universal constants) estimates of the expected values of the
mentioned quantities. 

The obtained formula may be also used in order to obtain exact formulas for
higher moments.

Let us comment on the organization of the paper.
In the next section we introduce some notation and prove the existence of moment-generating functions of truncated variation, upward truncated variation and downward truncated variation for any complex argument.
In the third section we calculate formula for the Laplace transform with respect to time parameter of the moment-generating function of upward truncated variation. In the fourth section we give examples of applications of the derived formula. 
In the last section we give possible interpretation of upward truncated variation in financial mathematics.

\section{Existence of moment-generating functions for any complex argument}

Let us start with some definitions and notation. The drawdown
and drawup processes of $W_{t}$ are defined respectively as 
\begin{eqnarray*}
DD_{s} &=&\sup_{0\leq t\leq s}W_{t}-W_{s}, \\
DU_{s} &=&W_{s}-\inf_{0\leq t\leq s}W_{t}.
\end{eqnarray*}
The times of drawdown of $c$ units and drawup of $c$ units are defined
respectively as 
\begin{eqnarray*}
T_{D}\left( c\right) &=&\inf \left\{ s\geq 0|DD_{s}=c\right\} , \\
T_{U}\left( c\right) &=&\inf \left\{ s\geq 0|DU_{s}=c\right\} .
\end{eqnarray*}

Further let $T_D^{\sup}(c)$ be the last instant when the maximum of $W_{t}$
on the interval $\left[ 0,T_D(c)\right] $ is attained and let $T_D^{\inf }(c)\leq T_D^{\sup}(c)$ be such that $W_{T_D^{\inf }(c)}=\inf_{0\leq s\leq T_D^{\sup}(c)}W_{s}.$

Let us fix $\alpha>0.$ We will prove the existence of moment-generating function of truncated variation, upward truncated variation and downward truncated variation for argument $\alpha$. Since the truncated variation and two other variables are non-negative, this will prove the existence of moment-generating function of those variables for any complex argument. 

\begin{proof}
Let $\delta >0$ be such a small number that 
\begin{equation*}
1-\mathbf{E}\exp \left( \alpha \sup_{0\leq t\leq T}W_{t}+\alpha c\right) \mathbf{P}\left( T_D(c)<\delta \right) >0.
\end{equation*}

By definition of $T_D(c)$ and $T_D^{\inf }(c)$ we have $W_{T_D^{\inf }(c)}>-c$
and hence, $W_{T_D^{\sup}(c)}-W_{T_D^{\inf }(c)}-c\leq W_{T_D^{\sup}(c)}.$ Let us fix $M>0.$ By Lemma 1 and Lemma 2 in \cite{L2008}, by independence of $W_{t}-W_{T_D(c)},t\geq T_D(c),$ and $T_D(c)$ (strong Markov property of Brownian motion) and by shift invariance property of truncated variation for stopping time $T_D(c)$ we have 
\begin{multline*}
\mathbf{E}\exp \left( \alpha TV_{\mu }^{c}\left[ 0,T\right] \wedge M\right)
\leq \mathbf{E}\exp \left( \alpha W_{T_D^{\sup}(c)}+\alpha c+\alpha TV_{\mu
}^{c}\left[ T_D(c),T\right] \wedge M\right) \\
\leq \mathbf{E}\exp \left( \alpha W_{T_D^{\sup}(c)}+\alpha c\right) \mathbf{E%
}\exp \left[ \alpha TV_{\mu }^{c}\left[ T_D(c),T\right] \wedge M;T_D(c)<\delta \right] \\
+\mathbf{E}\exp \left( \alpha W_{T_D^{\sup}(c)}+\alpha c\right) \mathbf{E}%
\exp \left[ \alpha TV_{\mu }^{c}\left[ T_D(c),T\right] \wedge M;T_D(c)\geq
\delta \right] \\
\leq \mathbf{E}\exp \left( \alpha W_{T_D^{\sup}(c)}+\alpha c\right) \mathbf{E%
}\exp \left[ \alpha TV_{\mu }^{c}\left[ T_D(c),T + T_D(c)\right] \wedge M;T_D(c)<\delta \right] \\
+\mathbf{E}\exp \left( \alpha W_{T_D^{\sup}(c)}+\alpha c\right) \mathbf{E}%
\exp \left[ \alpha TV_{\mu }^{c}\left[ T_D(c),T + T_D(c) - \delta \right] \wedge M;T_D(c)\geq \delta \right] \\
\leq \mathbf{E}\exp \left( \alpha \sup_{0\leq t\leq T}W_{t}+\alpha c\right) 
\mathbf{E}\exp \left( \alpha TV_{\mu }^{c}\left[ 0,T\right] \wedge M\right) 
\mathbf{P}\left( T_D(c)<\delta \right) \\
+\mathbf{E}\exp \left( \alpha \sup_{0\leq t\leq T}W_{t}+\alpha c\right) 
\mathbf{E}\exp \left( \alpha TV_{\mu }^{c}\left[0 ,T -\delta \right] \wedge
M\right) \mathbf{P}\left( T_D(c)\geq \delta \right) .
\end{multline*}
From the above we have 
\begin{multline*}
\mathbf{E}\exp \left( \alpha TV_{\mu }^{c}\left[ 0,T\right] \wedge M\right)
\\
\leq \frac{\mathbf{E}\exp \left( \alpha \sup_{0\leq t\leq T}W_{t}+\alpha
c\right) \mathbf{P}\left( T_D(c)\geq \delta \right) }{1-\mathbf{E}\exp \left( \alpha
\sup_{0\leq t\leq T}W_{t}+\alpha c\right) \mathbf{P}\left( T_D(c)<\delta \right) }%
\mathbf{E}\exp \left( \alpha TV_{\mu }^{c}\left[0 ,T -  \delta\right] \wedge
M\right) .
\end{multline*}
Similarly 
\begin{multline*}
\mathbf{E}\exp \left( \alpha TV_{\mu }^{c}\left[0,T -  \delta \right] \wedge
M\right)\\
\leq \frac{\mathbf{E}\exp \left( \alpha \sup_{0\leq t\leq T}W_{t}+\alpha
c\right) \mathbf{P}\left( T_D(c)\geq \delta \right) }{1-\mathbf{E}\exp \left( \alpha
\sup_{0\leq t\leq T}W_{t}+\alpha c\right) \mathbf{P}\left( T_D(c)<\delta \right) }%
\mathbf{E}\exp \left( \alpha TV_{\mu }^{c}\left[ 0 ,T-2\delta \right]
\wedge M\right) .
\end{multline*}
Iterating and putting together the above inequalities we finally obtain 
\begin{equation*}
\mathbf{E}\exp \left( \alpha TV_{\mu }^{c}\left[ 0,T\right] \wedge M\right)
\leq \left( \frac{\mathbf{E}\exp \left( \alpha \sup_{0\leq t\leq
T}W_{t}+\alpha c\right) \mathbf{P}\left( T_D(c)\geq \delta \right) }{1-\mathbf{E}%
\exp \left( \alpha \sup_{0\leq t\leq T}W_{t}+\alpha c\right) \mathbf{P}\left(
T_D(c)<\delta \right) }\right) ^{\left\lceil T/\delta \right\rceil }.
\end{equation*}
Letting $M\rightarrow \infty $ we get $\mathbf{E}\exp \left( \alpha TV_{\mu
}^{c}\left[ 0,T\right] \right) <+\infty .$

By (\ref{rel}) we obtain the finiteness of moment-generating functions 
of $UTV_{\mu }^{c}\left[ 0,T\right] $ and $DTV_{\mu }^{c}\left[ 0,T\right] $
as well.
\end{proof}

\section{Calculation of the Laplace transform of the moment-generating function}
\label{Calc}

Due to typographical reasons let us
introduce notation $\max \left\{ x,0\right\} =:(x)_{+}.$

The main difference between truncated variation and upward as well as downward truncated variation is such that for the latter quantities we have the following analog of Lemma 2 from \cite{L2008}, where instead of inequality we have equality.

\begin{lemma}
We have the following identities 
\begin{equation}
UTV_{\mu }^{c}\left[ 0,T\right] =\sup_{0\leq t<s\leq T_{D}(c)\wedge T}\left(
W_{s}-W_{t}-c\right) _{+}+UTV_{\mu }^{c}\left[ T_{D}(c)\wedge T,T\right] .  \label{n8}
\end{equation}
and
\begin{equation}
DTV_{\mu }^{c}\left[ 0,T\right] =\sup_{0\leq t<s\leq T_{U}(c)\wedge T}\left(
W_{t}-W_{s}-c\right) _{+}+DTV_{\mu }^{c}\left[ T_{U}(c)\wedge T,T\right] . 
\end{equation}
\end{lemma}

\begin{proof}
We will only prove the first formula (\ref{n8}), since the proof of the second one is identical. 

Let $0\leq t_{1}<s_{1}<t_{2}<s_{2}...<t_{n}<s_{n}\leq T$ be numbers\ from
the interval $\left[ 0,T\right] .$

We will prove that 
\begin{equation}
\sum_{i=1}^{n}\left( W_{s_{i}}-W_{t_{i}}-c\right) _{+}\leq \sup_{0\leq
t<s\leq T_{D}(c)\wedge T}\left( W_{s}-W_{t}-c\right) _{+}+UTV_{\mu }^{c}\left[
T_{D}(c)\wedge T,T\right] .  \label{n7}
\end{equation}
Let $n_{0}$ be the greatest number such that $s_{n_{0}}<T_{D}(c)$ and let us
assume that $n_{0}<n$ and $t_{n_{0}+1}<T_{D}(c).$

Let us consider several cases.

\begin{itemize}
\item  $W_{t_{n_{0}+1}}\geq W_{T_{D}(c)}.$ In this case 
\begin{equation*}
\left( W_{s_{n_{0}+1}}-W_{t_{n_{0}+1}}-c\right) _{+}\leq \left(
W_{s_{n_{0}+1}}-W_{T_{D}(c)}-c\right) _{+}.
\end{equation*}
and 
\begin{align}
\sum_{i=1}^{n}\left( W_{s_{i}}-W_{t_{i}}-c\right) _{+}& \leq
\sum_{i=1}^{n_{0}}\left( W_{s_{i}}-W_{t_{i}}-c\right) _{+}+\left(
W_{s_{n_{0}+1}}-W_{T_{D}(c)}-c\right) _{+}  \notag \\
& +\sum_{i=n_{0}+2}^{n}\left( W_{s_{i}}-W_{t_{i}}-c\right) _{+}.  \label{n4}
\end{align}

\item  $W_{t_{n_{0}+1}}<W_{T_{D}(c)}$ and $W_{s_{n_{0}+1}}\leq W_{T_{D}(c)^{\sup
}}.$ In this case $t_{n_{0}+1}<T_{D}^{\sup}(c)$ (since for $T_{D}^{\sup}(c) <t<T_{D}(c),$ $W_{t}>W_{T_{D}(c)}$) so 
\begin{equation*}
\left( W_{s_{n_{0}+1}}-W_{t_{n_{0}+1}}-c\right) _{+}\leq \left(
W_{T_{D}^{\sup }(c)}-W_{t_{n_{0}+1}}-c\right) _{+}
\end{equation*}
and 
\begin{align}
\sum_{i=1}^{n}\left( W_{s_{i}}-W_{t_{i}}-c\right) _{+}& \leq
\sum_{i=1}^{n_{0}}\left( W_{s_{i}}-W_{t_{i}}-c\right) _{+}+\left(
W_{T_{D}^{\sup }(c)}-W_{t_{n_{0}+1}}-c\right) _{+}  \notag \\
& +\sum_{i=n_{0}+2}^{n}\left( W_{s_{i}}-W_{t_{i}}-c\right) _{+}.  \label{n5}
\end{align}

\item  $W_{t_{n_{0}+1}}<W_{T_{D}(c)}$ and $W_{s_{n_{0}+1}}>W_{T_{D}^{\sup
}(c)}=W_{T_{D}(c)}+c.$ In this case 
\begin{multline*}
\left( W_{s_{n_{0}+1}}-W_{t_{n_{0}+1}}-c\right)
_{+}=W_{s_{n_{0}+1}}-W_{t_{n_{0}+1}}-c \\
=W_{T_{D}^{\sup }(c)}-W_{t_{n_{0}+1}}-c+W_{s_{n_{0}+1}}-W_{T_{D}^{\sup }(c)} \\
=W_{T_{D}^{\sup }(c)}-W_{t_{n_{0}+1}}-c+W_{s_{n_{0}+1}}-W_{T_{D}(c)}-c \\
=\left( W_{T_{D}^{\sup }(c)}-W_{t_{n_{0}+1}}-c\right) _{+}+\left(
W_{s_{n_{0}+1}}-W_{T_{D}(c)}-c\right) _{+}
\end{multline*}
and 
\begin{align}
\sum_{i=1}^{n}\left( W_{s_{i}}-W_{t_{i}}-c\right) _{+}& \leq
\sum_{i=1}^{n_{0}}\left( W_{s_{i}}-W_{t_{i}}-c\right) _{+}+\left(
W_{T_{D}^{\sup }(c)}-W_{t_{n_{0}+1}}-c\right) _{+}  \notag \\
& +\left( W_{s_{n_{0}+1}}-W_{T_{D}(c)}-c\right) _{+}+\sum_{i=n_{0}+2}^{n}\left(
W_{s_{i}}-W_{t_{i}}-c\right) _{+}.  \label{n6}
\end{align}
\end{itemize}

Thus for $t_{n_{0}+1}<T_{D}(c)$\ inequality (\ref{n4}), (\ref{n5}) or (\ref{n6}) holds and we may assume, adding in the case $t_{n_{0}+1}<T_{D}(c)$ new terms
in the partition and renaming the old ones, that 
\begin{eqnarray*}
0 &\leq &t_{1}<s_{1}<...<t_{n_{0}}<s_{n_{0}}\leq T_{D}(c), \\
T_{D}(c) &\leq &t_{n_{0}+1}<s_{n_{0}+1}<...<t_{n}<s_{n}\leq T.
\end{eqnarray*}

In order to prove (\ref{n7}) without loss of generality we may assume that
for any $1\leq i\leq n_{0},$ $\left( W_{s_{i}}-W_{t_{i}}-c\right) _{+}>0$
(otherwise we may omit the summand $\left( W_{s_{i}}-W_{t_{i}}-c\right) _{+}
$). From definition of $T_{D}(c)$ we have that for any $1\leq i\leq n_{0}-1$, $%
W_{s_{i}}-W_{t_{i+1}}<c,$ so 
\begin{align*}
& \left( W_{s_{i}}-W_{t_{i}}-c\right) _{+}+\left(
W_{s_{i+1}}-W_{t_{i+1}}-c\right) _{+} \\
& =W_{s_{i}}-W_{t_{i}}-c+W_{s_{i+1}}-W_{t_{i+1}}-c \\
& =W_{s_{i+1}}-W_{t_{i}}-c+\left( W_{s_{i}}-W_{t_{i+1}}-c\right)
<W_{s_{i+1}}-W_{t_{i}}-c.
\end{align*}
Iterating the above inequality, we obtain 
\begin{equation*}
\sum_{i=1}^{n_{0}}\left( W_{s_{i}}-W_{t_{i}}-c\right) _{+}\leq
W_{s_{n_{0}}}-W_{t_{1}}-c\leq \sup_{0\leq t<s\leq T_{D}(c)\wedge T}\left(
W_{s}-W_{t}-c\right) _{+}.
\end{equation*}
This, together with the obvious inequality 
\begin{equation*}
\sum_{i=n_{0}+1}^{n}\left( W_{s_{i}}-W_{t_{i}}-c\right) _{+}\leq \ UTV_{\mu
}^{c}\left[ T_{D}(c)\wedge T,T\right] 
\end{equation*}
proves (\ref{n7}). Taking supremum over all partitions $0\leq
t_{1}<s_{1}<t_{2}<s_{2}<...<t_{n}<s_{n}\leq T$ we finally get 
\begin{equation*}
UTV_{\mu }^{c}\left[ 0,T\right] \leq \sup_{0\leq t<s\leq T_{D}(c)\wedge T}\left(
W_{s}-W_{t}-c\right) _{+}+UTV_{\mu }^{c}\left[ T_{D}(c)\wedge T,T\right] .
\end{equation*}
Since the opposite inequality is obvious, we finally get (\ref{n8}).
\end{proof}

Now we are ready to state

\begin{lemma}
Let $\lambda $ be an arbitrary complex number and let 
\begin{equation*}
L(\lambda ,T):=\mathbf{E}\exp (\lambda UTV_{\mu }^{c}[0,T]),
\end{equation*}
$T>0,$ be a family of moment-generating functions of variables $UTV_{\mu }^{c}\left[
0,T\right] .$ This family satisfies the following integral equation 
\begin{multline}
L(\lambda ,T)=\int_{0}^{T}\int_{c}^{\infty }e^{\lambda \left( y-c\right)
}L(\lambda ,T-t)\mathbf{P}\left( T_{D}\left( c\right) \in dt,\sup_{0\leq
s\leq T_{D}\left( c\right) }DU_{s}\in dy\right) \\
+\int_{0}^{T}L(\lambda ,T-t)\mathbf{P}\left( T_{D}\left( c\right) \in
dt,\sup_{0\leq s\leq T_{D}\left( c\right) }DU_{s}<c\right) \\
+\int_{c}^{\infty }e^{\lambda \left( y-c\right) }\mathbf{P}\left(
T_{D}\left( c\right) >T,\sup_{0\leq s\leq T}DU_{s}\in dy\right)  \label{n1}
\\
+\mathbf{P}\left( T_{D}\left( c\right) >T,\sup_{0\leq s\leq
T}DU_{s}<c\right) .
\end{multline}
\end{lemma}

\begin{proof}
By Lemma 3 we have that for any $T>0$ 
\[
UTV_{\mu }^{c}\left[ 0,T\right] =\sup_{0\leq s\leq T_{D}\left( c\right)
\wedge T}\left( DU_{s}-c\right) _{+}+UTV_{\mu }^{c}\left[ T_{D}\left(
c\right) \wedge T,T\right].
\]

From dependence of $W_{t},t\in \left[ 0,T_{D}\left( c\right) \wedge T\right]
$ and $W_{t}-W_{T_{D}\left( c\right) \wedge T},t\in \left[ T_{D}\left(
c\right) \wedge T,T\right] ,$ only through $T_{D}\left( c\right) ,$ and by equality of distribution of $UTV_{\mu }^{c}\left[ T_{D}\left(
c\right) \wedge T,T\right]$ and $UTV_{\mu }^{c}\left[ 0,T-T_{D}\left(
c\right) \wedge T\right]$ we have

\begin{multline}
\mathbf{E}\exp \left( \mathbf{\lambda }UTV_{\mu }^{c}\left[ 0,T\right]
\right)  \\
=\mathbf{E}\exp \left( \lambda \sup_{0\leq s\leq T_{D}\left( c\right) \wedge
T}\left( DU_{s}-c\right) _{+}+\lambda UTV_{\mu }^{c}\left[ T_{D}\left(
c\right) \wedge T,T\right] \right)  \\
=\int_{0}^{\infty }\mathbf{E}\exp \left( \lambda \sup_{0\leq s\leq t\wedge
T}\left( DU_{s}-c\right) _{+}\right) \mathbf{E}\exp \left( \lambda UTV_{\mu
}^{c}\left[ 0,T-t\wedge T\right] \right) \mathbf{P}\left( T_{D}\left(
c\right) \in dt\right)  \\
=\int_{0}^{T}\int_{c}^{\infty }e^{\lambda \left( y-c\right) }\mathbf{E}\exp
\left( \lambda UTV_{\mu }^{c}\left[ 0,T-t\right] \right) \mathbf{P}\left(
T_{D}\left( c\right) \in dt,\sup_{0\leq s\leq T_{D}\left( c\right)
}DU_{s}\in dy\right)  \\
+\int_{0}^{T}\mathbf{E}\exp \left( \lambda UTV_{\mu }^{c}\left[ 0,T-t\right]
\right) \mathbf{P}\left( T_{D}\left( c\right) \in dt,\sup_{0\leq s\leq
T_{D}\left( c\right) }DU_{s}<c\right)  \\
+\int_{c}^{\infty }e^{\lambda \left( y-c\right) }\mathbf{P}\left(
T_{D}\left( c\right) >T,\sup_{0\leq s\leq T}DU_{s}\in dy\right)  \\
+\mathbf{P}\left( T_{D}\left( c\right) >T,\sup_{0\leq s\leq
T}DU_{s}<c\right) .  \notag
\end{multline}

In the third line of the calculations above we have used iterated
expectation, strong Markov property and the shift invariance of upward
truncated variation for stopping time $T_{D}(c).$

\end{proof}

Hadjiliadis and Zhang in their recent paper (\cite{HZ2009}) calculated for 
$a,b>0$ the densities 
\begin{equation*}
p\left( t;a,b\right) dt = \mathbf{P}\left( T_{D}\left( a\right) \in
dt,T_{U}\left( b\right) >t\right)
\end{equation*}
and 
\begin{equation*}
q\left( t;a,b\right) dt = \mathbf{P}\left( T_{U}\left( a\right) \in
dt,T_{D}\left( b\right) >t\right).
\end{equation*}
Using these densities we are able to write equation (\ref{n1}) in more
elegant form. Indeed, we have

\begin{lemma}
The family $L(\lambda ,T)$ satisfies the following integral equation 
\begin{multline}
L\left( \lambda ,T\right) =\int_{0}^{T}L\left( \lambda ,T-t\right) \left\{
p\left( t;c,c\right) +\int_{c}^{\infty }e^{\lambda \left( y-c\right) }\frac{%
\partial p\left( t;c,y\right) }{\partial y}dy\right\} dt \\
-\int_{0}^{T}\mathbf{P}\left( T_{D}\left( c\right) >T-t\right) \left\{
q\left( t;c,c\right) +\int_{c}^{\infty }e^{\lambda \left( y-c\right) }\frac{%
\partial q\left( t;y,c\right) }{\partial y}dy\right\} dt \\
+\mathbf{P}\left( T_{D}\left( c\right) >T\right) .  \label{r5}
\end{multline}
\end{lemma}

\begin{proof}
We have 
\begin{multline}
\mathbf{P}\left( T_{D}\left( c\right) \in dt,\sup_{0\leq s\leq T_{D}\left(
c\right) }DU_{s}\in dy\right) \\
=\mathbf{P}\left( T_{D}\left( c\right) \in dt,T_{U}\left( y+dy\right)
>t\right) -\mathbf{P}\left( T_{D}\left( c\right) \in dt,T_{U}\left( y\right)
>t\right)  \label{n2} \\
=\frac{\partial p\left( t;c,y\right) }{\partial y}dydt
\end{multline}
and 
\begin{eqnarray}
\mathbf{P}\left( T_{D}\left( c\right) \in dt,\sup_{0\leq s\leq T_{D}\left(
c\right) }DU_{s}<c\right) &=&\mathbf{P}\left( T_{D}\left( c\right) \in
dt,T_{U}\left( c\right) >t\right) \notag \\
&=&p\left( t;c,c\right) dt. \label{nn2}
\end{eqnarray}

In order to express $\mathbf{P}\left(
T_{D}\left( c\right) >T,\sup_{0\leq s\leq T}DU_{s}\in dy\right)$ with $p(t;a,b)$ and $q(t;a,b)$ let us notice that for $y>0$
\begin{align}
& \mathbf{P}\left( T_{D}\left( c\right) >T,\sup_{0\leq s\leq T}DU_{s}\geq
y\right) \notag \\
& =\int_{0}^{T}\mathbf{P}\left( T_{U}\left( y\right) \in dt,T_{D}\left(
c\right) >T\right)\notag \\
& =\int_{0}^{T}\mathbf{P}\left( T_{U}\left( y\right) \in dt,T_{D}\left(
c\right) >t\right) \mathbf{P}\left( T_{D}\left( c\right) >T-t\right) \notag \\
& =\int_{0}^{T}q\left( t;y,c\right) \mathbf{P}\left( T_{D}\left( c\right)
>T-t\right) dt \label{nn4}
\end{align}
The equality 
\begin{align*}
& \mathbf{P}\left( T_{U}\left( y\right) \in dt,T_{D}\left( c\right) >T\right) \\
& = \mathbf{P}\left( T_{U}\left( y\right) \in dt,T_{D}\left( c\right) >t\right) 
\mathbf{P}\left( T_{D}\left( c\right) >T-t\right)
\end{align*}
holds since the event $\left\{ T_{U}\left( y\right) \in dt\right\} $ also
means that the process $W_{t}$ reaches a new maximum at the moment $t.$ Now
for $y>0$ we calculate 
\begin{multline}
\mathbf{P}\left( T_{D}\left( c\right) >T,\sup_{0\leq s\leq T}DU_{s}\in
dy\right) \\
=\mathbf{P}\left( T_{D}\left( c\right) >T,\sup_{0\leq s\leq T}DU_{s}\geq
y\right) -\mathbf{P}\left( T_{D}\left( c\right) >T,\sup_{0\leq s\leq
T}DU_{s}\geq y+dy\right) \\
=\int_{0}^{T}\left\{ q\left( t;y,c\right) -q\left( t;y+dy,c\right) \right\} 
\mathbf{P}\left( T_{D}\left( c\right) >T-t\right) dt \\
=-\int_{0}^{T}\frac{\partial q\left( t;y,c\right) }{\partial y}\mathbf{P}%
\left( T_{D}\left( c\right) >T-t\right) dtdy.  \label{n3}
\end{multline}
Using similar reasoning, by (\ref{nn4}) we also have 
\begin{multline}
\mathbf{P}\left( T_{D}\left( c\right) >T,\sup_{0\leq s\leq T}DU_{s}<c\right)
=\mathbf{P}\left( T_{D}\left( c\right) >T,T_{U}\left( c\right) >T\right) \\
=\mathbf{P}\left( T_{D}\left( c\right) >T\right) -\int_{0}^{T}q\left(
t;c,c\right) \mathbf{P}\left( T_{D}\left( c\right) >T-t\right) dt. \label{nn3}
\end{multline}
Thus, from (\ref{n1}), (\ref{n2}), (\ref{nn2}), (\ref{n3}) and (\ref{nn3}) we obtain the integral equation (\ref{r5}) satisfied by the family of moment-generating functions of upward truncated variation: 
\begin{eqnarray}
L\left( \lambda ,T\right) &=&\int_{0}^{T}\int_{c}^{\infty }e^{\lambda \left(
y-c\right) }L\left( \lambda ,T-t\right) \frac{\partial p\left( t;c,y\right) 
}{\partial y}dydt  \notag \\
&&+\int_{0}^{T}L\left( \lambda ,T-t\right) p\left( t;c,c\right) dt  \notag \\
&&-\int_{c}^{\infty }e^{\lambda \left( y-c\right) }\int_{0}^{T}\frac{%
\partial q\left( t;y,c\right) }{\partial y}\mathbf{P}\left( T_{D}\left(
c\right) >T-t\right) dtdy  \notag \\
&&+\quad\mathbf{P}\left( T_{D}\left( c\right) >T\right) -\int_{0}^{T}q\left(
t;c,c\right) \mathbf{P}\left( T_{D}\left( c\right) >T-t\right) dt  \notag \\
&=&\int_{0}^{T}L\left( \lambda ,T-t\right) \left\{ p\left( t;c,c\right)
+\int_{c}^{\infty }e^{\lambda \left( y-c\right) }\frac{\partial p\left(
t;c,y\right) }{\partial y}dy\right\} dt  \notag \\
&&-\int_{0}^{T}\mathbf{P}\left( T_{D}\left( c\right) >T-t\right) \left\{
q\left( t;c,c\right) +\int_{c}^{\infty }e^{\lambda \left( y-c\right) }\frac{%
\partial q\left( t;y,c\right) }{\partial y}dy\right\} dt  \notag \\
&&+\quad\mathbf{P}\left( T_{D}\left( c\right) >T\right) .  \notag
\end{eqnarray}
\end{proof}

In order to shorten notation let introduce new functions of parameters $t$
and $\lambda$ 
\begin{gather*}
p\left(\lambda, t\right) :=p\left( t;c,c\right) +\int_{0}^{\infty
}e^{\lambda y}\frac{\partial p\left( t;c,y+c\right) }{\partial y}dy, \\
q\left(\lambda, t\right) :=q\left( t;c,c\right) +\int_{0}^{\infty
}e^{\lambda y}\frac{\partial q\left( t;y+c,c\right) }{\partial y}dy
\end{gather*}
and for such pairs of complex numbers $(\lambda, \nu)$ that the integral $%
\int_{0}^{\infty }e^{\nu t}L\left( \lambda ,t\right)dt$ exists, let us
define 
\begin{gather*}
M\left(\lambda, \nu \right) :=\int_{0}^{\infty }e^{\nu t}L\left( \lambda
,t\right) dt, \\
T\left( \nu \right) :=\int_{0}^{\infty }e^{\nu t}\mathbf{P}\left(
T_{D}\left( c\right) >t\right) dt.
\end{gather*}
By (\ref{r5}) we have 
\begin{eqnarray*}
M\left(\lambda, \nu \right) &=&\int_{0}^{\infty }e^{\nu \tau
}L\left(\lambda, \tau \right) d\tau =\int_{0}^{\infty }e^{\nu \tau
}\int_{0}^{\tau }L\left(\lambda, \tau -t\right)p\left(\lambda, t\right)
dtd\tau \\
&&-\int_{0}^{\infty }e^{\nu \tau }\int_{0}^{\tau }\mathbf{P}\left(
T_{D}\left( c\right) >\tau -t\right) q\left(\lambda, t\right) dtd\tau
+T\left( \nu \right) \\
&=&\int_{0}^{\infty }e^{\nu t}p\left(\lambda, t\right) \int_{t}^{\infty
}e^{\nu \left( \tau -t\right) }L\left(\lambda, \tau -t\right) d\tau dt \\
&&-\int_{0}^{\infty }e^{\nu t}q\left(\lambda, t\right) \int_{t}^{\infty
}e^{\nu \left( \tau -t\right) }\mathbf{P}\left( T_{D}\left( c\right) >\tau
-t\right) d\tau dt+T\left( \nu \right) \\
&=&M\left(\lambda, \nu \right) \int_{0}^{\infty }e^{\nu t}p\left(\lambda,
t\right) dt-T\left( \nu \right) \int_{0}^{\infty }e^{\nu t}q\left(\lambda,
t\right) dt+T\left( \nu \right) .
\end{eqnarray*}
Thus we obtained a formula for the Laplace transform with respect to $T$ of
the moment-generating function of $UTV_{\mu }^{c}\left[ 0,T\right] :$ 
\begin{equation}
M\left(\lambda, \nu \right) =T\left( \nu \right) \frac{1-\int_{0}^{\infty
}e^{\nu t}q\left(\lambda, t\right) dt}{1-\int_{0}^{\infty }e^{\nu
t}p\left(\lambda, t\right) dt}.  \label{r7}
\end{equation}

Using results of \cite{HZ2009} and \cite{Taylor1975} we are able to compute $%
M\left( \lambda ,\nu \right) $ more directly. We have

\begin{theorem}
For $\nu $ with negative real part and any complex $\lambda$ the following formula holds 
\begin{equation}
M\left( \lambda ,\nu \right) =-\frac{1}{\nu }-\frac{\lambda e^{\mu c}}{\nu
^{2}}\frac{\mu \sinh \left( cU_{\mu }\left( \nu \right) \right) -U_{\mu
}\left( \nu \right) \cosh \left( cU_{\mu }\left( \nu \right) \right) }{\frac{%
\lambda U_{\mu }\left( \nu \right) }{\nu }+\sinh \left( 2cU_{\mu }\left( \nu
\right) \right) -2\frac{\lambda +\mu }{U_{\mu }\left( \nu \right) }\sinh
^{2}\left( cU_{\mu }\left( \nu \right) \right) }, \label{th}
\end{equation}
where $U_{\mu }\left( \nu \right) =\sqrt{\mu ^{2}-2\nu }.$
\end{theorem}

\begin{proof}
Integrating by parts, we obtain 
\begin{eqnarray}
T\left( \nu \right)  &=&\int_{0}^{\infty }e^{\nu t}\mathbf{P}\left(
T_{D}\left( c\right) >t\right) dt  \notag \\
&=&\frac{e^{\nu t}}{\nu }\mathbf{P}\left( T_{D}\left( c\right) >t\right)
|_{t=0}^{t=\infty }-\int_{0}^{\infty }\frac{e^{\nu t}}{\nu }\frac{d}{dt}%
\mathbf{P}\left( T_{D}\left( c\right) >t\right) dt  \notag \\
&=&-\frac{1}{\nu }+\frac{1}{\nu }\mathbf{E}e^{\nu T_{D}\left( c\right) }.
\label{r1}
\end{eqnarray}
Similarly, we have 
\begin{eqnarray*}
p\left( \lambda ,t\right)  &=&p\left( t;c,c\right) +\int_{0}^{\infty
}e^{\lambda y}\frac{\partial p\left( t;c,y+c\right) }{\partial y}dy \\
&=&-\lambda \int_{0}^{\infty }e^{\lambda y}p\left( t;c,y+c\right) dy,
\end{eqnarray*}
hence 
\begin{eqnarray}
\int_{0}^{\infty }e^{\nu t}p\left( \lambda ,t\right) dt &=&-\lambda
\int_{0}^{\infty }\int_{0}^{\infty }e^{\lambda y}p\left( t;c,y+c\right) dydt
\notag \\
&=&-\lambda \int_{0}^{\infty }e^{\lambda y}\int_{0}^{\infty }e^{\nu t}%
\mathbf{P}\left( T_{D}\left( c\right) \in dt,T_{U}\left( y+c\right)
>t\right) dy  \notag \\
&=&-\lambda \int_{0}^{\infty }e^{\lambda y}\mathbf{E}e^{\nu T_{D}\left(
c\right) }I_{\left\{ T_{U}\left( y+c\right) >T_{D}\left( c\right) \right\}
}dy.  \label{r2}
\end{eqnarray}
Using notation from \cite{HZ2009}, page 11, we have 
\begin{equation*}
\mathbf{E}e^{\nu T_{D}\left( c\right) }I_{\left\{ T_{U}\left( y+c\right)
>T_{D}\left( c\right) \right\} }=\left( 1-L_{0}^{-W}\left( -\nu ;c\right)
e^{T_{-\mu ,1}\left( -\nu ,c\right) y}\right) \mathbf{E}e^{\nu T_{D}\left(
c\right) }
\end{equation*}
thus 
\begin{align*}
& \int_{0}^{\infty }e^{\lambda y}\mathbf{E}e^{\nu T_{D}\left( c\right)
}I_{\left\{ T_{U}\left( y+c\right) >T_{D}\left( c\right) \right\} }dy \\
& =\left( \int_{0}^{\infty }e^{\lambda y}\left[ 1-L_{0}^{-W}\left( -\nu
,c\right) \exp \left( T_{-\mu ,1}\left( -\nu ,a\right) y\right) \right]
dy\right) \mathbf{E}e^{\nu T_{D}\left( c\right) } \\
& =\left( \frac{L_{0}^{-W}\left( -\nu ,c\right) }{T_{-\mu ,1}\left( -\nu
,c\right) +\lambda }-\frac{1}{\lambda }\right) \mathbf{E}e^{\nu T_{D}\left(
c\right) }
\end{align*}
and finally from (\ref{r2}) we obtain 
\begin{equation}
\int_{0}^{\infty }e^{\nu t}p\left( \lambda ,t\right) dt=\left( 1-\lambda 
\frac{L_{0}^{-W}\left( -\nu ,c\right) }{T_{-\mu ,1}\left( -\nu ,c\right)
+\lambda }\right) \mathbf{E}e^{\nu T_{D}\left( c\right) }.  \label{r4}
\end{equation}
Similarly 
\begin{eqnarray*}
q\left( \lambda ,t\right)  &=&q\left( t;c,c\right) +\int_{0}^{\infty
}e^{\lambda y}\frac{\partial q\left( t;y+c,c\right) }{\partial y}dy \\
&=&-\lambda \int_{0}^{\infty }e^{\lambda y}q\left( t;y+c,c\right) dy,
\end{eqnarray*}
hence 
\begin{eqnarray*}
\int_{0}^{\infty }e^{\nu t}q\left( \lambda ,t\right) dt &=&-\lambda
\int_{0}^{\infty }e^{\nu t}\int_{0}^{\infty }e^{\lambda y}q\left(
t;y+c,c\right) dydt \\
&=&-\lambda \int_{0}^{\infty }e^{\lambda y}\mathbf{E}e^{\nu T_{U}\left(
y+c\right) }I_{\left\{ T_{U}\left( y+c\right) <T_{D}\left( c\right) \right\}
}dy.
\end{eqnarray*}
Again, by results of \cite{HZ2009} and using symmetry of standard Brownian motion, we have 
\begin{equation*}
\mathbf{E}e^{\nu T_{U}\left( y+c\right) }I_{\left\{ T_{U}\left( y+c\right)
<T_{D}\left( c\right) \right\} }=L_{0}^{-W}\left( -\nu ;c\right) e^{T_{-\mu
,1}\left( -\nu ;c\right) y},
\end{equation*}
and finally we get 
\begin{eqnarray}
\int_{0}^{\infty }e^{\nu t}q\left( \lambda ,t\right) dt &=&-\lambda
\int_{0}^{\infty }e^{\lambda y}\mathbf{E}e^{\nu T_{U}\left( y+c\right)
}I_{\left\{ T_{U}\left( y+c\right) <T_{D}\left( c\right) \right\} }dy  \notag
\\
&=&-\lambda \int_{0}^{\infty }e^{\lambda y}L_{0}^{-W}\left( -\nu ,c\right)
e^{T_{-\mu ,1}\left( -\nu ,c\right) y}dy  \notag \\
&=&\lambda \frac{L_{0}^{-W}\left( -\nu ,c\right) }{T_{-\mu ,1}\left( -\nu
,c\right) +\lambda }.  \label{r3}
\end{eqnarray}
Finally from (\ref{r7}), (\ref{r1}), (\ref{r4}) and (\ref{r3}) we obtain 
\begin{eqnarray}
M\left( \lambda ,\nu \right)  &=&\left( -\frac{1}{\nu }+\frac{1}{\nu }%
\mathbf{E}e^{\nu T_{D}\left( c\right) }\right) \frac{1-\lambda \frac{%
L_{0}^{-W}\left( -\nu ,c\right) }{T_{-\mu ,1}\left( -\nu ,c\right) +\lambda }%
}{1-\left( 1-\lambda \frac{L_{0}^{-W}\left( -\nu ,c\right) }{T_{-\mu
,1}\left( -\nu ,c\right) +\lambda }\right) \mathbf{E}e^{\nu T_{D}\left(
c\right) }}  \notag \\
&=&-\frac{1}{\nu }\frac{\left( 1-\lambda \frac{L_{0}^{-W}\left( -\nu
,c\right) }{T_{-\mu ,1}\left( -\nu ,c\right) +\lambda }\right) \left( 1-%
\mathbf{E}e^{\nu T_{D}\left( c\right) }\right) }{1-\left( 1-\lambda \frac{%
L_{0}^{-W}\left( -\nu ,c\right) }{T_{-\mu ,1}\left( -\nu ,c\right) +\lambda }%
\right) \mathbf{E}e^{\nu T_{D}\left( c\right) }}.  \label{r6}
\end{eqnarray}

It is possible to express the obtained formula for $M\left( \lambda ,\nu
\right) $ with the elementary functions. We have (cf. \cite{HZ2009} and \cite{Taylor1975}): 
\begin{gather*}
L_{0}^{-W}\left( -\nu ,c\right) =\frac{U_{\mu }\left( \nu \right) }{-2\nu }%
\left\{ \frac{e^{\mu c}\left( U_{\mu }\left( \nu \right) \coth \left(
cU_{\mu }\left( \nu \right) \right) -\mu \right) }{\sinh \left( cU_{\mu
}\left( \nu \right) \right) }-\frac{U_{\mu }\left( \nu \right) }{\sinh
^{2}\left( cU_{\mu }\left( \nu \right) \right) }\right\} , \\
\mathbf{E}e^{\nu T_{D}\left( c\right) }=\frac{U_{\mu }\left( \nu \right)
e^{-\mu c}}{U_{\mu }\cosh \left( cU_{\mu }\left( \nu \right) \right) -\mu
\sinh \left( cU_{\mu }\left( \nu \right) \right) }
\end{gather*}
and 
\begin{equation*}
T_{-\mu ,1}\left( -\nu ,c\right) =\mu -U_{\mu }\left( \nu \right) \coth
\left( cU_{\mu }\left( \nu \right) \right) ,
\end{equation*}
where 
\begin{equation*}
U_{\mu }\left( \nu \right) =\sqrt{\mu ^{2}-2\nu }.
\end{equation*}
Substituting the above formulas in (\ref{r6}) we obtain (\ref{th}).
\end{proof}

\section{Examples of applications}

\label{Appl}The direct application of the derived formula may be calculation
of the moment-generating function $L\left( \lambda ,T\right) $ (with the use of the
inverse Laplace transform formula). However, we will start with simpler
formulae.

\subsection{Exact formula for the expected value of $UTV_{\mu }^{c}%
\left[ 0,T\right] $.}

Using formula 
\begin{equation*}
\mathbf{E}UTV_{\mu }^{c}\left[ 0,T\right] =\lim_{\lambda \rightarrow 0}%
\frac{1}{\lambda }\left( L\left( \lambda ,T\right) -1\right)
\end{equation*}
we obtain 
\begin{align}
& \lim_{\lambda \rightarrow 0}\frac{1}{\lambda }\left( M\left( \nu
,\lambda \right) -M\left( \nu ,0\right) \right)   \notag \\
& =\lim_{\lambda \rightarrow 0}\frac{1}{\lambda }\left( \int_{0}^{\infty
}e^{\nu t}L\left( \lambda ,t\right) dt-\int_{0}^{\infty }e^{\nu t}L\left(
0,t\right) dt\right)   \notag \\
& =\int_{0}^{\infty }e^{\nu t}\lim_{\lambda \rightarrow 0}\frac{1}
{\lambda }\left[ L\left( \lambda ,t\right) -1\right] dt  \notag \\
& =\int_{0}^{\infty }e^{\nu t}\mathbf{E}UTV_{\mu }^{c}\left[ 0,t\right] dt.
\label{r8}
\end{align}
On the other hand, from (\ref{r6}) we have 
\begin{align}
& \lim_{\lambda \rightarrow 0}\frac{1}{\lambda }\left( M\left( \nu
,\lambda \right) -M\left( \nu ,0\right) \right)  \notag \\
& =\lim_{\lambda \rightarrow 0}\frac{1}{\lambda }\left( -\frac{1}{\nu }%
\frac{\left( 1-\lambda \frac{L_{0}^{-W}\left( -\nu ,c\right) }{T_{-\mu
,1}\left( -\nu ,c\right) +\lambda }\right) \left( 1-\mathbf{E}e^{\nu
T_{D}\left( c\right) }\right) }{1-\left( 1-\lambda \frac{L_{0}^{-W}\left(
-\nu ,c\right) }{T_{-\mu ,1}\left( -\nu ,c\right) +\lambda }\right) \mathbf{E%
}e^{\nu T_{D}\left( c\right) }}+\frac{1}{\nu }\right)  \notag \\
& =\frac{L_{0}^{-W}\left( -\nu ,c\right) }{\nu T_{-\mu ,1}\left( -\nu
,c\right) \left( 1-\mathbf{E}e^{\nu T_{D}\left( c\right) }\right) },
\label{r9}
\end{align}
which, by (\ref{r8}) and after substituting in (\ref{r9}) the formulas for $L_{0}^{-W}\left(
-\nu ,c\right),$ $\mathbf{E}e^{\nu T_{D}\left( c\right) }$ and $T_{-\mu
,1}\left( -\nu ,c\right) $ yields 
\begin{equation}
\int_{0}^{\infty }e^{\nu t}\mathbf{E}UTV_{\mu }^{c}\left[ 0,t\right] dt =\frac{e^{\mu c}\sqrt{\mu
^{2}-2\nu }}{2\nu ^{2}\sinh \left( c\sqrt{\mu ^{2}-2\nu }\right) }.
\label{r10}
\end{equation}

Inverting the formula (\ref{r10}) we are able to obtain exact formula for
the expected value of $\mathbf{E}UTV_{\mu }^{c}\left[ 0,T\right] .$ Let $%
\mathcal{L}_{s}^{-1}(g)$ denote inverse of the Laplace transform of the
function $g(s)=\int_{0}^{\infty }e^{-s\cdot t}f\left( t\right) dt,$ i.e. the
function $f\left( t\right) .$ We have 
\begin{equation}
\mathcal{L}_{s}^{-1}\left( s^{-2}\right) =t  \label{r15}
\end{equation}
and, by the last formula on page 641 of \cite{BS2002} 
\begin{equation*}
\mathcal{L}_{s}^{-1}\left( \frac{\sqrt{2s}}{\sinh \left( c\sqrt{2s}\right) }%
\right) =\frac{\sqrt{2}}{\sqrt{\pi }t^{5/2}}\sum_{k=0}^{\infty }\left(
\left( 2k+1\right) ^{2}c^{2}-t\right) e^{-\left( 2k+1\right)
^{2}c^{2}/\left( 2t\right) }.
\end{equation*}
Hence, by properties of Laplace transform 
\begin{align}
& \mathcal{L}_{s}^{-1}\left( \frac{\sqrt{2s+\mu ^{2}}}{\sinh \left( c\sqrt{%
2s+\mu ^{2}}\right) }\right)   \notag \\
& =\frac{\sqrt{2}}{\sqrt{\pi }t^{5/2}}e^{-\mu ^{2}t}\sum_{k=0}^{\infty
}\left( \left( 2k+1\right) ^{2}c^{2}-t\right) e^{-\left( 2k+1\right)
^{2}c^{2}/\left( 2t\right) }.  \label{r16}
\end{align}
Finally, by (\ref{r15}), (\ref{r16}) and Borel convolution theorem for the
Laplace transform, we obtain 
\begin{align}
& \mathbf{E}UTV_{\mu }^{c}\left[ 0,T\right]   \notag \\
& =\mathcal{L}_{s}^{-1}\left( \frac{e^{\mu c}\sqrt{2s+\mu ^{2}}}{2s^{2}\sinh
\left( c\sqrt{2s+\mu ^{2}}\right) }\right)   \notag \\
& =\frac{e^{\mu c}}{\sqrt{2\pi }}\int_{0}^{T}\left( T-t\right) \frac{e^{-\mu
^{2}t}}{t^{5/2}}\sum_{k=0}^{\infty }\left( \left( 2k+1\right)
^{2}c^{2}-t\right) e^{-\left( 2k+1\right) ^{2}c^{2}/\left( 2t\right) }dt 
\notag \\
& =\frac{e^{\mu c}}{\sqrt{2\pi }}\sum_{k=0}^{\infty }\int_{0}^{T}\left(
T-t\right) \frac{\left( 2k+1\right) ^{2}c^{2}-t}{t^{5/2}}e^{-\mu
^{2}t-\left( 2k+1\right) ^{2}c^{2}/\left( 2t\right) }dt.  \label{r18}
\end{align}

\subsection{Estimation of the expected value of $UTV_{\protect\mu }^{c}\left[
0,T\right] .$}

In \cite{L2008} we obtained a formula
for function $F\left( \mu ,c,T\right) ,$ such that \ 
\begin{equation}
\mathbf{E}TV_{\mu }^{c}\left[ 0,T\right] \sim F\left(  \left| \mu \right| ,c,T\right),
\label{r20}
\end{equation}
where relation ''$\sim $'' means that the ratio $\mathbf{E}TV_{\mu }^{c}%
\left[ 0,T\right] /F\left(\left| \mu \right|,c,T\right) $ is separated from $0$ and infinity by universal constants, which do not depend on $\mu, c, T$.

On the other hand, we see that the exact formula (\ref{r18}) for $\mathbf{E}UTV_{\mu }^{c}\left[
0,T\right] $ may be stated in the form 
\begin{equation}
\mathbf{E}UTV_{\mu }^{c}\left[ 0,T\right] =e^{\mu c}G\left( \left| \mu
\right| ,c,T\right) ,  \label{r19}
\end{equation}
where 
\begin{equation*}
G\left( \left| \mu \right| ,c,T\right) =\frac{1}{\sqrt{2\pi }}%
\sum_{k=0}^{\infty }\int_{0}^{T}\left( T-t\right) \frac{\left( 2k+1\right)
^{2}c^{2}-t}{t^{5/2}}e^{-\mu ^{2}t-\left( 2k+1\right) ^{2}c^{2}/\left(
2t\right) }dt
\end{equation*}
does not depend on the sign of $\mu .$ Using (\ref{rel}), (\ref{r19}) and the fact that $DTV_{\mu }^{c}\left[ 0,T\right] $ has the same distribution as $UTV_{-\mu}^{c}\left[ 0,T\right] $ we see that 
\begin{eqnarray}
\mathbf{E}TV_{\mu }^{c}\left[ 0,T\right] &\sim &\mathbf{E}UTV_{\mu }^{c}%
\left[ 0,T\right] +\mathbf{E}DTV_{\mu }^{c}\left[ 0,T\right]  \notag \\
&=&\mathbf{E}UTV_{\mu }^{c}\left[ 0,T\right] +\mathbf{E}UTV_{-\mu }^{c}\left[
0,T\right]  \notag \\
&\sim &e^{\left| \mu \right| c}G\left( \left| \mu \right| ,c,T\right).
\label{r21}
\end{eqnarray}

Comparing (\ref{r20}) and (\ref{r21}) we see that 
\begin{equation*}
G\left( \left| \mu \right| ,c,T\right) \sim e^{-\left| \mu \right| c}F\left(
\left| \mu \right| ,c,T\right)
\end{equation*}
and finally we get estimates up to universal constants for $\mathbf{E}%
UTV_{\mu }^{c}\left[ 0,T\right] :$ 
\begin{align*}
& \mathbf{E}UTV_{\mu }^{c}\left[ 0,T\right] \sim e^{\mu c-\left| \mu
\right| c}F\left( \left| \mu \right| ,c,T\right) \\
& = e^{\mu c-\left| \mu \right| c}\left\{ 
\begin{array}{cl}
T/c+ \left| \mu \right| T & \text{if }\sqrt{T}\geq\chi(c, \mu); \\ 
2\sqrt{T}+ \left| \mu \right| T-c & \text{if } \sqrt{T} \in \left(c- \left| \mu \right|  T , \chi(c, \mu) \right); \\ 
T^{3/2}\frac{\exp(-\left( c- \left| \mu \right| T\right) ^{2}/(2T))} {\left( c- \left| \mu \right| T\right) ^{2}} & \text{if }\sqrt{T} \leq c- \left| \mu \right|  T \text{,}
\end{array}
\right.
\end{align*}
where $\chi \left( c, \mu \right) =\sqrt{\frac{e^{2\mu \left|
c\right| }-2\mu \left| c\right| -1}{2\mu ^{2}}}=c\sqrt{1+\frac{2}{3} \left| \mu \right| c+...} \geq c.$

\subsection{Laplace transform of the second moment\textbf{\ }of $UTV_{%
\protect\mu }^{c}\left[ 0,T\right] $}

Similarly as (\ref{r18}) we may obtain a formula for the Laplace transform
of the second moment of $UTV_{\mu }^{c}\left[ 0,T\right] :$ 
\begin{align}
& \int_{0}^{\infty }e^{\nu t}\mathbf{E}\left( UTV_{\mu }^{c}\left[ 0,t\right]
\right) ^{2}dt=\left[ \frac{\partial ^{2}}{\partial \lambda ^{2}}M\left( \nu
,\lambda \right) \right] _{\lambda =0}  \notag \\
& =-\frac{2L_{0}^{-W}\left( -\nu ,c\right) \left( 1-\mathbf{E}e^{\nu
T_{D}\left( c\right) }+L_{0}^{-W}\left( -\nu ,c\right) \mathbf{E}e^{\nu
T_{D}\left( c\right) }\right) }{\nu \left( T_{-\mu ,1}\left( -\nu ,c\right)
\left( 1-\mathbf{E}e^{\nu T_{D}\left( c\right) }\right) \right) ^{2}}.
\label{r17}
\end{align}

After substituting in formula (\ref{r17}) the formulas for $L_{0}^{-W}\left(
-\nu ,c\right) ,\mathbf{E}e^{\nu T_{D}\left( c\right) }$ and $T_{-\mu
,1}\left( -\nu ,c\right) ,$ it simplifies to 
\begin{align*}
& \int_{0}^{\infty }e^{\nu t}\mathbf{E}\left( UTV_{\mu }^{c}\left[ 0,t\right]
\right) ^{2}dt \\
& =-\frac{e^{\mu c}U_{\mu }\left( \nu \right) \left[ U_{\mu }^{2}\left( \nu
\right) +\nu \left( 1-\cosh \left( 2cU_{\mu }\left( \nu \right) \right)
\right) \right] }{2\nu ^{3}\left[ U_{\mu }\left( \nu \right) \cosh \left(
cU_{\mu }\left( \nu \right) \right) -\mu \sinh \left( cU_{\mu }\left( \nu
\right) \right) \right] \sinh ^{2}\left( cU_{\mu }\left( \nu \right) \right) 
}.
\end{align*}

\begin{remark}
Using formulas from \cite{BS2002} (page 642) it is possible to invert the
above formula and obtain expression for $\mathbf{E}\left( UTV_{\mu }^{c}%
\left[ 0,t\right] \right) ^{2}$ in terms of parabolic cylinder functions.
\end{remark}

\section{Interpretation of upward truncated variation in financial mathematics}

As it was mentioned earlier, upward truncated variation appears naturally in the expression for {\bf the least upper bound} for the rate of return from any trading of a financial asset, dynamics of which follows geometric Brownian motion, in the presence of flat commission. Similar result was proved in \cite{L2008} for truncated variation, however, truncated variation is not the least upper bound.

Indeed, similarly as in  \cite{L2008}, let us assume that the dynamics of the prices $P_t$ of some financial asset (e.g. stock) is the following $P_{t}=\exp \left(\mu t+\sigma B_{t} \right)$.
We are interested in the maximal possible profit coming from trading
this single instrument during time interval $\left[ 0,T\right].$ We buy the instrument at the moments $0\leq t_{1}<...<t_{n}<T$ and sell it at the moments $s_{1}<...<s_{n}\leq T,$ such that $t_{1}<s_{1}<t_{2}<s_{2}<...<t_{n}<s_{n},$ in order to obtain the maximal possible profit. Furthermore we assume that for every transaction we have to pay a flat commission and $\gamma $ is the ratio of the transaction value paid for the commission. 

The maximal possible rate of return from our strategy reads as (cf. \cite{L2008})
\[
\sup_{n}\sup_{0\leq t_{1}<s_{1}<...<t_{n}<s_{n}\leq T}\frac{%
P_{s_{1}}}{P_{t_{1}}}\frac{1-\gamma }{1+\gamma }...\frac{%
P_{s_{n}}}{P_{t_{n}}}\frac{1-\gamma }{1+\gamma } - 1. 
\]

Let $M_n$ be the set of all partitions $$\pi = \left\{ 0\leq t_{1}<s_{1}<...<t_{n}<s_{n}\leq T\right\}.$$ 
To see that $\exp \left(
\sigma UTV_{\mu /\sigma }^{c/\sigma }\left[ 0,T\right] \right)-1$ with $c=\ln 
\frac{1+\gamma }{1-\gamma }$\ is the least upper bound for maximal possible rate of return let us substitute 
\begin{equation*}
\begin{split}
\sup_{n} & \sup_{M_{n}}\prod_{i=1}^{n}\left\{ \frac{P_{s_{i}}}{P_{t_{i}}}\frac{1-\gamma }
{1+\gamma }\right\} 
 = \sup_{n}\sup_{M_{n}}\prod_{i=1}^{n}\left\{ \frac{\exp \left( \mu s_{i}+\sigma
B_{s_{i}}\right) }{\exp \left( \mu t_{i}+\sigma B_{t_{i}}\right) 
}e^{-c}\right\} \\
&=\sup_{n}\sup_{M_{n}}\exp \left( \sigma \sum_{i=1}^{n}\left\{ \left( \frac{\mu}{\sigma }s_{i}+B_{s_{i}}\right) - \left( \frac{\mu }{\sigma }
t_{i}+B_{t_{i}}\right) -\frac{c}{\sigma }\right\} \right) \\
&= \exp \left( \sigma \sup_{n}\sup_{M_{n}}\sum_{i=1}^{n}\left\{
\left( \frac{\mu }{\sigma }s_{i}+B_{s_{i}}\right) - \left( \frac{\mu 
}{\sigma }t_{i}+B_{t_{i}}\right) -\frac{c}{\sigma }\right\} \right)
\\
&= \exp \left( \sigma UTV_{\mu /\sigma }^{c/\sigma }\left[ 0,T\right]
\right) .
\end{split}
\end{equation*}
This gives the claimed bound.

\end{document}